\documentclass [12 pt, leqno] {amsart}

\usepackage {amssymb}
\usepackage {amsthm}
\usepackage {amscd}

\usepackage[all]{xy}
\usepackage {hyperref}

\begin {document}

\newtheorem{proposition}[equation]{Proposition}
\newtheorem{lemma}[equation]{Lemma}
\newtheorem{corollary}[equation]{Corollary}
\newtheorem{thm}[equation]{Theorem} 
\newtheorem*{thm*}{Theorem}

\theoremstyle{definition}

\theoremstyle{plain}

\numberwithin{equation}{section}

\newcommand{\eq}[2]{\begin{equation}\label{#1}  #2 \end{equation}}

\newcommand{\ml}[2]{\begin{multline}\label{#1}  #2 \end{multline}}

\newcommand{\nnal}[1]{\begin{align*} #1 \end{align*}}
\newcommand{\nneq}[1]{\begin{equation} \nonumber #1 \end{equation}}

\newcommand{\nnml}[1]{\begin{multline}\nonumber #1 \end{multline}}

\newcommand{\arir}{\ar@{^{(}->}}
\newcommand{\aril}{\ar@{_{(}->}}

\newcommand{\are}{\ar@{>>}}

\newcommand{\xr}[1] {\xrightarrow{#1}}




\newcommand{\codim}{{\rm codim}}

\newcommand{\Spec} {{\rm Spec}}

\newcommand{\Hom} {{\rm Hom}}

\newcommand{\rk} {{\rm rk}}


\newcommand{\Z} {\mathbb{Z}}
\newcommand{\Q} {\mathbb{Q}}
\newcommand{\C} {\mathbb{C}}
\newcommand{\G} {\mathbb{G}}

\newcommand{\A}{\mathbb{A}}
\renewcommand{\P}{\mathbb{P}}





\newcommand{\Hm} [3] {H^{#1}(#2,\Z(#3))}
\newcommand{\M}[1]{M_{gm}(#1)}

\title{Motivic cohomology of the complement of hyperplane arrangements}
\author{Andre Chatzistamatiou }
\thanks{This work has been supported by the DFG Leibniz Program}

\begin{abstract}
We give a presentation of the motivic cohomology ring of the complement of a hyperplane arrangement considered as algebra over the motivic
cohomology of the ground field.   
\end{abstract}
\maketitle

\section*{Introduction}
Let $K$ be a field and $U\subset \A_K^N$ the complement of a finite union of hyperplanes. For $K=\C$ the cohomology ring
$H^*(U^{an},\Z_{U^{an}})$ of the constant sheaf $\Z_{U^{an}}$ is isomorphic to the subalgebra of meromorphic forms generated 
by the  logarithmic forms $\frac{1}{2\pi i} \frac{df}{f}$  (\cite{Arnold},\cite{Brieskorn}). 
Motivic cohomology $H(U):=\oplus_{p,q} H^p(U,\Z(q))$ as defined by Voevodsky \cite{Voevodsky} is defined over an arbitrary 
field and there is a natural isomorphism $\G_{m}(U) \cong \Hm{1}{U}{1}$. 
The aim of this article is to give a description of the motivic cohomology ring $H(U)$ for a perfect ground field $K$.
We describe the module structure over $H(K)$ and the ring structure. 

Concerning the module structure, we show that $H(U)$ is a free $H(K)$ module
(Corollary \ref{freierModul}) and give a combinatorial
description for the rank (Corollary \ref{Rang}). 
For the ring structure, we consider $H(U)$ as an algebra over $H(K)$. We prove that $H(U)$ is generated by the units
$\G_m(U)$ and give the relations. We have to introduce some notation to describe the relations and to state the theorem. 
Let $H(K)\{\G_m(U)\}$ be the free, bigraded commutative algebra over $H(K)$ generated by the abelian group of units of $U$ 
(in degree $(1,1)$). If $f\in \G_m(U)$ then we denote by $(f)$ the corresponding element in $H(K)\{\G_m(U)\}$ and for 
$\lambda \in K^{\times}$  we denote by $[\lambda]\in \Hm{1}{K}{1}\subset H(K)$ the associated cohomology class. 

\begin{thm*}[Theorem \ref{thm}] 
Let $K$ be a perfect field. The morphism 
\nneq
{
 H(K)\{\G_m(U)\}/I \xr{} H(U),
}
defined by mapping $(f)$ to the class $[f]$ of $f$ in $\Hm{1}{U}{1}$, is an isomorphism of $H(K)$ algebras, where  
the ideal $I$ is generated by the elements:
\begin{small}
\nnal
{
&(1) \quad (f)-[f], \quad \text{if $f\in K^{\times} \subset \mathbb{G}_m(U)$,} \\
&(2) \quad (f_1)\cdot (f_2)\cdots (f_t), \quad \text{if $f_i\in  \mathbb{G}_m(U), i=1,\dots,t,$ s.t. $\sum_{i=1}^t{f_k}=1$,} \\
&(3) \quad (f)^2+[-1]\cdot (f), \quad \text{if $f\in \mathbb{G}_m(U)$.} 
} 
\end{small}  
\end{thm*}
\quad \\

Relation (1) is trivial, relation (2) and (3) are well-known in the Milnor K-theory of the function field $K(U)$ 
(\cite{Milnor}, Lemma 1.2, 1.3). \\
\quad \\
\emph{Acknowledgements.} This paper covers parts of my PhD thesis with supervisor H\'el\`ene Esnault. I am very grateful 
to her for the excellent guidance and the support. I would like to thank Marc Levine, Stefan Kukulies and Kay R\"ulling   
for helpful discussions and the referee for helpful comments.

\section{Generating classes}

\subsection{} 
Let $K$ be a perfect field. We work with the triangulated category $DM^{eff}_{gm}(K)$ of effective
geometrical motives over $K$, defined by Voevodsky (\cite{Voevodsky},
Definition 2.1.1). 
We denote by $M_{gm}$ the functor $Sm/K\xr{} DM^{eff}_{gm}(K)$ that maps a smooth scheme  to the 
corresponding motive. 
For a smooth scheme $X$ of finite type over $K$ motivic
cohomology is defined as 
\nneq
{
 H^p(X,\Z(q)):=\Hom_{DM^{eff}_{gm}}(M_{gm}(X),\Z(q)[p]).
}

Let  $U\subset \A_K^N$ be the complement of a finite union of
hyperplanes. 

\begin{proposition}\label{prop:reinpp}
There is an isomorphism 
\nneq
{
 M_{gm}(U)\cong \oplus_{i\in I} \Z(n_i)[n_i]
}
for some finite index set $I$ and and integers $n_i\geq 0$.
\begin{proof}
Let $Y_1,\dots,Y_r$ be the hyperplanes in the complement of $U$, i.e. $U=\A^N-\cup_{i=1}^r Y_i$. We proceed by induction on $r$. 
If $U=\A^N_K$ then $M_{gm}(U)=\Z$.  If $r\geq 1$ then the assertion is true for $U':=\A^N-\cup_{i=2}^r Y_i$ and 
$Y_1\cap U'=Y_1-\cup_{i=2}^r(Y_1\cap  Y_i)$ by induction. Consider the Gysin triangle 
\small 
\eq{Spaltungmotivisch}
 {M_{gm}(U) \xr{} M_{gm}(U') \xr{\phi} M_{gm}(Y_1\cap U')(1)[2] \xr{+1}. 
} \normalsize
 We have $\phi=0$ since 
\eq{claimcan}
{
 \Hom_{DM^{eff}_{gm}}(\Z(n)[n],\Z(m)[m+1])=0 \quad \text{for every $n,m\geq 0$,}
}
which is proved in Lemma (\ref{vanishing}) below.
 It follows that the sequence (\ref{Spaltungmotivisch}) is split and there is a non-canonical isomorphism 
\eq{unnatuerlicherSpalt}
{
M_{gm}(U)\cong M_{gm}(U') \oplus M_{gm}(Y_1\cap U')(1)[1],
}  
which proves the assertion.
\end{proof}
\end{proposition}

\begin{lemma} \label{vanishing}
For every $n,m\geq 0$ the identity (\ref{claimcan}) holds.
\begin{proof}
By the Cancellation Theorem (\cite{VoevodskyCancellation}, Corollary 4.10) we may reduce to 
\nnal
{
 &\Hom_{DM^{eff}_{gm}}(\Z(n)[n],\Z[1])=0 , \quad n\geq 0 \\
 &\Hom_{DM^{eff}_{gm}}(\Z,\Z(n)[n+1])=0 , \quad n\geq 0. 
}
The first group is a direct summand of 
\nneq
{
\Hom_{DM^{eff}_{gm}}(M_{gm}(\G_m^{\times n}),\Z[1])= H^1_{{\rm
 Zar}}(\G_m^{\times n},\Z) = 0;
}
the second group is a direct summand of
$H^1_{Zar}(K,C_{*}(\G_m^{\times n}))$, which is trivial because $C_{*}(\G_m^{\times n})$ is
concentrated in non-positive degrees. 
\end{proof}
\end{lemma}

The splitting of the Gysin sequence (\ref{Spaltungmotivisch}) yields the \emph{split} short 
exact sequence
\small 
\eq{SpaltungKohomologie}
{
 0\xr{} \Hm{p}{U'}{q} \xr{} \Hm{p}{U}{q} \xr{} \Hm{p-1}{Y_1\cap U'}{q-1}  \xr{} 0.
} \normalsize 

Using induction again and the isomorphism (\ref{unnatuerlicherSpalt}) we have:
\begin{corollary}\label{freierModul}
The motivic cohomology $H(U)$ is a  finitely generated free module over the motivic cohomology of the ground field $H(K)$.
\end{corollary}
 
\subsection{} A unit $f\in \G_m(U)$ gives a cohomology class $[f]\in \Hm{1}{U}{1}$ defined by the morphism 
\nneq
{
 M_{gm}(U)\xr{f} M_{gm}(\G_m) \xr{\cong} \Z \oplus \Z(1)[1] \xr{{\rm proj}} \Z(1)[1].  
}
The decomposition $M_{gm}(\G_m)\cong \Z\oplus \Z(1)[1]$ is constructed from the splitting 
$\Z \xr{=} M_{gm}(K) \xr{\text{inclusion at $1$}} \M{\G_m}$. 
By (\cite{Voevodsky}, Corollary 3.5.3) the map $f\mapsto [f]$ is an isomorphism of abelian groups. 

\begin{proposition}\label{proposition:Erzeuger} 
The cohomology ring $H(U)$ is generated by the classes of the units in $U$ as an algebra over $H(K)$. 
\begin{proof}
We use the short exact sequence (\ref{SpaltungKohomologie})
\small 
\nneq 
{
 0\xr{} \oplus_{p,q} \Hm{p}{U'}{q} \xr{\alpha} \oplus_{p,q} \Hm{p}{U}{q} \xr{\beta}
 \oplus_{p,q} \Hm{p-1}{Y_1\cap U'}{q-1}  \xr{} 0
} \normalsize
of $H(K)$ modules. The map $\alpha$ is the restriction to $U$, hence it is a ring homomorphism which maps 
$[f]$ to $[f \mid_U]$ for $f\in \G_m(U')$.  For the map $\beta$, let $t\in \A^1(U')$  be a function with vanishing locus
$Y_1\cap U'$ and denote by $\iota:H(U') \xr{} H(Y_1\cap U')$ the restriction morphism. I claim that 
\eq{claim}
{
\beta( [t]\cdot \alpha(z) ) = \iota(z)
}
holds for any $z$. In $DM^{eff}_{gm}$ we need to show that the composition 
\small
\nnal
{
 M_{gm}(Y_1\cap U')(1)[1]\xr{\text{Gysin}} \M{U} \xr{(incl,t)} \M{ U' \times \G_m } \\
                                                 \xr{id\otimes {\rm proj}} \M{U'}(1)[1]  
} \normalsize
is equal to  $M_{gm}({\rm inclusion})(1)[1]$. This follows from the morphism of Gysin triangles:
\small
$$
\begin{CD}
\M{Y_1\cap U'}(1)[1] @>{+1}>>\M{U} @>>> \M{U'} @>>>  \\
@VV{incl(1)[1]}V @VV{(incl,t)}V @VV{(id,t)}V  \\
\M{U'}(1)[1] @>{+1}>> \M{U'\times \G_m} @>>> \M{U'\times \A^1} @>>>
\end{CD}
$$\normalsize
and the commutative diagram
$$
\xymatrix
{
\M{U'}(1)[1] \ar[r]^{\text{Gysin}} \ar[dr]_{=}
&
\M{U'\times \G_m} \ar[d]^{id_{U'}\otimes {\rm proj} }
\\
&
\M{U'}(1)[1].
}
$$

In particular, if $f_1,\dots,f_s\in \G_m(U')$ then 
\nneq
{
\beta([t]\cdot [f_1] \cdots [f_s])= [f_1\mid_{Y_1\cap U'}]\cdots [f_s\mid_{Y_1\cap U'}].
}
Since $\G_m(U')\xr{} \G_m(Y_1\cap U')$ is surjective, we conclude by induction on the number of removed hyperplanes 
that the algebra generated by $\G_m(U)$ maps via $\beta$ surjectively onto $H(Y_1\cap U')$ and contains $H(U')$, thus equals $H(U)$. 
\end{proof}
\end{proposition}

\section{Relations}

\subsection{} 
Let $U$ be a smooth $K$-scheme. The purpose of this section is to show that the following elements in the 
motivic cohomology of $U$ are trivial:
\eq{Relation1}
{
[f_1]\cdot [f_2]\cdots [f_t]=0, \quad \text{if $f_i\in  \mathbb{G}_m(U), i=1\dots t$,}\; \text{such that $\sum_{i=1}^t{f_i}=1$,} 
}
and 
\eq{Relation2}
{
[f]^2+[-1]\cdot [f]=0, \quad \text{if $f\in \mathbb{G}_m(U)$.} 
} 
To prove the identities, we reduce the general case to the case $U=\Spec(K)$; here we use the comparison theorem 
of motivic cohomology and Milnor $K$-theory by Suslin and Voevodsky (see \cite{SuslinVoevodsky})
\eq{SVM}
{
H^n(K,\Z(n))\cong K^M_n(K), \quad \text{for any $n\geq 0$.}
}  

It is sufficient to prove (\ref{Relation1}) for the scheme $\Delta^t-\cup_{i=1}^t \{x_i=0\}$, where  
$\Delta^t=\{\sum_i x_i=1\} \subset \A^{t}$, and $f_1=x_1,\dots,f_t=x_t$, because the general case follows
by pullback. Similarly, it is enough to show (\ref{Relation2}) for $\G_m$. It is somewhat easier to prove the identity
\eq{Relation3}
{
 R(f_1,\dots,f_t)=0 , \quad \text{if $f_i\in  \mathbb{G}_m(U), i=1\dots t$,}\;  \text{such that $\sum_{i=1}^t f_i=0$},
}
where $R(f_1,\dots,f_t):=$ 
\nnml
{
  \sum_{i=1}^t (-1)^{i} [f_1]\cdots  \widehat{[f_i]} \cdots  [f_t] 
 + 
 \sum_{i<j} [-1]\cdot
 [f_1]\cdots  \widehat{[f_i]} \cdots \widehat{[f_j]}  \cdots  [f_t] 
 + \\
 \sum_{i<j<k} [-1]^2 \cdot 
 [f_1] \cdots \widehat{[f_i]} \cdots \widehat{[f_j]} \cdots  \widehat{[f_k]}  \cdots [f_t] 
 + \dots. 
}

\begin{lemma} \label{EquivalenzRelationen}
The identities (\ref{Relation1}) together with  (\ref{Relation2}) are equivalent to the
identity (\ref{Relation3}).
\begin{proof}
First, assume (\ref{Relation1}) and (\ref{Relation2}), let $f_1,\dots,f_t$ be
units and $\sum_i f_i=0$. Then $\sum_{i=1}^{t-1} \frac{f_i}{-f_t}=1$ and 
\eq{produkt}
{
 \left([f_1]-[f_t]+[-1]\right)\cdots \left([f_{t-1}]-[f_t]+[-1]\right)=0
}
by (\ref{Relation1}). 
Using skew-commutativity: \small$[f][g]+[g][f]=0$, \normalsize for all units $f,g$,  
we see with the help of \small$[f_t]^2-[-1]\cdot [f_t]=0$ \normalsize 
(\ref{Relation2}) and \small$[-1]+[-1]=0$  \normalsize that (\ref{produkt}) is equal to $(-1)^{t}R(f_1,\dots,f_t)$. 

Now we assume (\ref{Relation3}). Since \small{$R(f_1,\dots,f_{t-1},-1)=(-1)^{t}[f_1]\cdots [f_{t-1}]$} we get (\ref{Relation1}). 
In order to show (\ref{Relation2}) we may reduce to $U=\G_m$. 
If $K\not= \mathbb{F}_2$ then there exist units $\lambda_1,\lambda_2,\lambda_3 \in K^{\times}$ such that $\lambda_1+\lambda_2+\lambda_3=0$; 
and we have 
\nneq
{
 R(\lambda_1 f,\lambda_2 f, \lambda_3f)=[f]^2+[-1]\cdot[f] + R(\lambda_1,\lambda_2,\lambda_3).
}
For the case $K=\mathbb{F}_2$, we use $M_{gm}(\G_m)=\Z\oplus \Z(1)[1]$ to see 
$H^2(\G_m,\Z(2))=H^2(\mathbb{F}_2,\Z(2)) \oplus H^1(\mathbb{F}_2,\Z(1))$,
which is zero by the isomorphism with Milnor $K$-theory (\ref{SVM}). 
\end{proof} 
\end{lemma}

\subsection{}
Let $U\subset \A^N_K$ be the complement of a finite union of hyperplanes $Y_1,\dots,Y_r$. 
We define $U_j:=\A^N_K-\cup_{\substack{ i=1  \\ i\not =j}}^r Y_i$. 

The next lemma will serve as a criterion for a cohomology class to be trivial. 

\begin{lemma}\label{Kriterium}
Any morphism $\phi: \M{U}\xr{} T$ in $DM^{eff}_{gm}$ such that 
\nneq
{
 \M{Y_j\cap U_j}(1)[1]\xr{\text{Gysin}} \M{U} \xr{\phi} T  
}
is trivial for every $j=1,\dots,r$, factors through $\M{K}$, i.e. there is a morphism $\psi: \M{K}\xr{} T$ and a commutative diagram 
$$
\xymatrix
{
\M{U} \ar[r]^-{\phi} \ar[d] 
&
T. 
\\
\M{K} \ar[ur]_{\psi}
&
}
$$
\begin{proof}
We prove by induction on $r$. The case $U=\A^N$ is obvious. 

By assumption and 
with the help of the Gysin triangle (\ref{Spaltungmotivisch}) the morphism
$\phi$ factors as $\phi: \M{U} \xr{\M{incl}} \M{U_1} \xr{\phi_1} T$ for some $\phi_1$. For every $j>1$ the diagram 
$$
\xymatrix
{
\M{U_j \cap Y_j}(1)[1] \ar[r]^-{\text{Gysin}} \ar[dd]^{\M{incl}(1)[1]}
&
\M{U} \ar[dr]^{\phi} \ar[dd]
& 
\\
&
&
T
\\
\M{(U_j\cup U_1) \cap Y_j}(1)[1] \ar[r]^-{\text{Gysin}} 
&
\M{U_1} \ar[ur]_{\phi_1}
&
}
$$
is commutative. 
Since $\M{(U_j\cup U_1) \cap Y_j}$ is isomorphic (by the morphism $\M{incl}$) to a direct summand of $\M{U_j \cap Y_j}$ 
(\ref{unnatuerlicherSpalt}) it follows from 
$\phi_1 \circ \text{Gysin} \circ \M{incl}(1)[1]= \phi \circ \text{Gysin} = 0$ that $\phi_1 \circ \text{Gysin}=0$.
Now, we can apply induction to $\phi_1$. 
\end{proof}
\end{lemma}

\begin{proposition} \label{proposition:Relationen} 
  Let $U$ be a smooth $K$-scheme of finite type and
  $f_1,\dots, f_t\in \G_m(U)$ units such that $\sum_{i=1}^t f_i = 0$. Then 
  $$R(f_1,\dots,f_t)=0.$$
 \begin{proof}
  We may assume $U:=H-\cup_{i=1}^t \{x_i=0\}$, $H\subset \A^t$ the hyperplane
  $\sum_{i=1}^t x_i=0$, and $f_i=x_i$ for every $i$, where $x_1,\dots,x_t$ are
  the coordinates of $\A^t$.    

  Denote by \small $\beta_j:\Hm{t-1}{U}{t-1}\xr{} \Hm{t-2}{Y_j\cap U_j}{t-2}$\normalsize,
  for every \small $j=1,\dots, t$\normalsize, the morphism from the Gysin sequence; here $Y_j=\{x_j=0\}$
  and $U_j=H-\cup_{i \not= j}\{x_i=0\}$.
  Formula (\ref{claim}) immediately implies 
\nneq
{
\beta_j(R(x_1,\dots,x_t))=(-1)^{j} R(x_1\mid_{Y_j\cap U_j},\dots, \widehat{x_j},\dots, x_t\mid_{Y_j\cap U_j}),
}  
the righthand side being zero by induction on $t$. Lemma (\ref{Kriterium})
  implies that \small $R(x_1,\dots,x_t)$ \normalsize  is the pullback of some class in 
  \small$\Hm{t-1}{K}{t-1}$\normalsize. If $K\not=\mathbb{F}_2$ then $U$ has a $K$ rational
  point. After pullback to this point we have to prove:
  \small $R(\lambda_1,\dots,\lambda_t)=0$ \normalsize in \small $\Hm{t-1}{K}{t-1}$ \normalsize for $\lambda_i\in K^{\times}$
  with $\sum_i \lambda_i=0$. 
 We may use the isomorphism (\ref{SVM}) to work with Milnor $K$-theory. In Milnor $K$-theory the formulas
  (\ref{Relation1}), (\ref{Relation2}) are well-known (\cite{Milnor}, Lemma 1.2, 1.3) and they yield  
 \small $R(\lambda_1,\dots,\lambda_t)=0$ \normalsize as in the proof of Lemma (\ref{EquivalenzRelationen}). 
\end{proof}  
\end{proposition}  

With the help of  Lemma (\ref{EquivalenzRelationen}) we have the following  corollary. 

\begin{corollary}\label{cor:Relationen}
Let $U$ be a smooth $K$-scheme.
\begin{enumerate}
\item If $f_1,\dots, f_t\in \G_m(U)$ are units such that $\sum_{i=1}^t f_i = 1$ then 
      $[f_1]\cdots [f_t]=0$. 
\item For every $f\in \G_{m}(U)$ we have $[f]^2+[-1]\cdot [f]=0$.
\end{enumerate}
\end{corollary}

\section{The cohomology ring}

\subsection{} 
Recall that $U\subset \A^N_K$ is the complement of a finite union of hyperplanes $Y_1,\dots,Y_r$.

We let $H(K)\{\G_m(U)\}$ be the free, bigraded commutative algebra over $H(K)$
generated by the abelian group of units of $U$ in degree $(1,1)$. For $f\in
\G_m(U)$ we denote by $(f)$ the corresponding element in $H(K)\{\G_m(U)\}$. 

We denote by $I_U\subset H(K)\{\G_m(U)\}$ the ideal generated by the following elements:
\begin{small}
\begin{align}
 (f)-&[f], \quad \text{if $f\in K^{\times}\subset \mathbb{G}_m(U)$,} \label{Gl0}  
\\ 
 (f_1)\cdot (f_2)\cdots (f_t), \quad &\text{if $f_i\in  \mathbb{G}_m(U),
  i=1,\dots,t,$ such that $\sum_{i=1}^t{f_k}=1$,}  \label{Gl1}
\\
 (f)^2+&[-1]\cdot (f), \quad \text{if $f\in \mathbb{G}_m(U)$.} \label{Gl2}
\end{align}
\end{small}

 As above we define $\tilde{R}(f_1,\dots,f_t)\in H(K)\{\G_m(U)\}$ to be the element 
 \begin{small}
 \nnml
 {
  \sum_{i=1}^t (-1)^{i} (f_1)\cdots  \widehat{(f_i)} \cdots  (f_t) 
 + 
 \sum_{i<j} [-1]\cdot
 (f_1)\cdots  \widehat{(f_i)} \cdots \widehat{(f_j)}  \cdots  (f_t) 
 + \\
 \sum_{i<j<k} [-1]^2 \cdot 
 (f_1) \cdots \widehat{(f_i)} \cdots \widehat{(f_j)} \cdots  \widehat{(f_k)}  \cdots (f_t) 
 + \dots
 }
 \end{small}

By the same arguments as in the proof of Lemma (\ref{EquivalenzRelationen}) we
can replace (\ref{Gl1}) in the definition of $I_U$ by 
\eq{Gl2R}
{
 \tilde{R}(f_1,\dots,f_t), \quad  \text{if $f_i\in  \mathbb{G}_m(U), i=1,\dots,t,$ s.t. $\sum_{i=1}^t{f_k}=0$.} 
}
For every hyperplane $Y_i$ we choose a polynomial $\phi_i$ of degree $1$,
which defines $Y_i$. To prove the main theorem it will be useful to work with
the ideal $I'_U\subset I_U$ generated by elements of the form (\ref{Gl0}),
(\ref{Gl2}) and elements as in (\ref{Gl2R}) with $f_j=\lambda_j\phi_{i_j}$ or $f_j=\lambda_j$  for
every $j$ and some $\lambda_j\in K^{\times}$ and some index $i_j$.
 
Proposition (\ref{proposition:Erzeuger}) and Corollary (\ref{cor:Relationen}) imply that we have surjective morphisms 
\nneq
{
 H(K)\{\G_m(U)\}/I'_U \xr{} H(K)\{\G_m(U)\}/I_U \xr{} H(U), 
} 
defined by mapping $(f)$ to $[f]$.

\begin{thm}\label{thm}
Let $K$ be a perfect field. The morphism of $H(K)$ algebras
\nneq
{
 H(K)\{\G_m(U)\}/I_U \xr{} H(U),
}
is an isomorphism.
\begin{proof}
Obviously it is sufficient to prove the assertion for $I'_U$.
The proof is done by induction on the number of removed hyperplanes.

Let  $U':=\A^N_K - \bigcup_{i\geq 2} Y_i$, then by induction we have
\nnal
{
 H(K)\{\G_m(U')\}/I'_{U'} &\xr{\cong} H(U'),\\
 H(K)\{\G_m(Y_1\cap U')\}/I'_{Y_1\cap U'} &\xr{\cong} H(Y_1\cap U').
}   
The restriction map $\G_m(U')\xr{} \G_m(Y_1\cap U')$ induces
\nneq
{
\iota:H(K)\{\G_m(U')\}/I'_{U'} \xr{} H(K)\{\G_m(Y_1\cap U')\}/I'_{Y_1\cap U'},
} 
which is the pullback morphism 
in motivic cohomology for the inclusion $Y_1\cap U'\xr{} U'$. Furthermore,
the restriction $\G_m(U')\xr{} \G_m(U)$ induces 
\nneq
{
 \widetilde{\alpha} : H(K)\{\G_m(U')\}/I'_{U'} \xr{} H(K)\{\G_m(U)\}/I'_U,
}
and we define $\widetilde{\beta}$ to be the composition
\nneq
{
H(K)\{\G_m(U)\}/I'_U\xr{} \oplus_{p,q} \Hm{p}{U}{q} \xr{} \oplus_{p,q} \Hm{p-1}{Y_1\cap U'}{q-1},
} 
where the second arrow comes from Gysin sequence. To see how
$\widetilde{\beta}$ maps, we may write every element $x\in H(K)\{\G_m(U)\}$
as $$x=(\phi_1)\cdot \widetilde{\alpha}(x_1)+ \widetilde{\alpha}(x_2).$$ This
can be done by using relation (\ref{Gl2}). Then formula (\ref{claim}) yields
$\widetilde{\beta}(x)=\iota(x_1)$.

Obviously we have a commutative diagram
\tiny
$$
\xymatrix
{
H(K)\{\G_m(U')\}/I'_{U'} \ar[r]^-{\widetilde{\alpha}} \ar[d]^{\cong}
&
H(K)\{\G_m(U)\}/I'_{U} \ar[r]^-{\widetilde{\beta}} \ar[d]
&
H(K)\{\G_m(Y_1\cap U')\}/I'_{Y_1\cap U'} \ar[d]^{\cong}
\\
\oplus_{p,q} \Hm{q}{U'}{p}  \ar[r]
&
\oplus_{p,q} \Hm{q}{U}{p} \ar[r]
&
\oplus_{p,q} \Hm{q-1}{Y_1\cap U'}{p-1}
}
$$ \normalsize
and since the sequence (\ref{SpaltungKohomologie}) is exact we have to prove
that the first row is exact, i.e. if $\iota(x_1)=0$ then 
$(\phi_1)\cdot \widetilde{\alpha}(x_1)$ is contained in the image of $\widetilde{\alpha}$. 
Then 
\nneq
{
H(K)\{\G_m(U)\}/I'_U \xr{} H(U) 
}
is an isomorphism. 

First we prove that the kernel of $\iota$ is generated by the elements
\small 
\ml{iotaKern}
{
 \tilde{R}(f_1,\dots,f_t), \; \text{with}\; f_j=\lambda_j \phi_{i_j}, i_j>1,\;
 \text{or}\; f_j=\lambda_j, \\  \text{such that}\;
 \sum_j (f_j\mid_{Y_1\cap U'})=0. 
}\normalsize
We denote by $J$ the ideal generated by elements of the form (\ref{iotaKern})
in $H(K)\{\G_m(U')\}$.
The image of $I'_{U'}+J$ in $H(K)\{\G_m(Y_1\cap U')\}$ is equal to $I'_{Y_1\cap U'}$. In order to see this
we note that an element in $I'_{Y_1\cap U'}$ of the form (\ref{Gl2}) lifts to an element in $I'_{U'}$ 
(since the restriction $\G_m(U')\xr{} \G_m(Y_1\cap U')$ is surjective), and an element in $I'_{Y_1\cap U'}$ 
of the form (\ref{Gl2R}) with $f_j=\lambda_j\cdot \phi_{i_j}\mspace{-8.0mu}\mid_{Y_1}$ or $f_j=\lambda_j$ 
and $\lambda_j\in K^{\times}$ for every $j$ lifts to an element (\ref{iotaKern}) in $J$.

Thus it is sufficient to prove
\nneq
{
\text{ker}\left(\G_m(U')\xr{}\G_m(Y_1\cap U') \right) \subset I'_{U'}+J.
}
It is easy to see that the kernel is generated by elements of the form 
\begin{enumerate}
\item $\lambda \cdot \frac{\phi_i}{\phi_j}$ with $i,j$ such that $Y_1\cap
  Y_i=Y_1\cap Y_j$; $\lambda=\frac{\phi_j\mid_{Y_1}}{\phi_i\mid_{Y_1}}$,
\item $\lambda \phi_i$ with $i$ such that $Y_i\cap Y_1=\emptyset$; $\lambda=\frac{1}{\phi_i\mid_{Y_1}}$.  
\end{enumerate}
Since \small{$\left(\frac{\lambda \phi_i}{\phi_j}\right)=\tilde{R}(\lambda\phi_i,-\phi_j)\in J$} and
\small{$(\lambda\phi_i)=\tilde{R}(\lambda\phi_i,-1)\in J$}  the kernel of $\iota$ is $J$ as claimed.

Let \small{$x_1=\tilde{R}(f_1,\dots,f_t)$} be as in (\ref{iotaKern}). 
We see that $\sum_{j} f_j= -\mu\cdot \phi_1$ 
for suitable $\mu\in K$ since every $f_j$ is a polynomial of degree $\leq 1$ and the restriction $\sum_{j} f_j\mspace{-8.0mu}\mid_{Y_1}$
to $Y_1$ is vanishing. 
The case $\mu=0$ is trivial, thus we may assume $\mu\not=0$. The following calculation in $H(K)\{\G_m(U)\}/I'_U$ completes the proof:
\begin{small}
\nnal
{
 (\phi_1)\cdot \widetilde{\alpha}(x_1)&= (\mu \phi_1)\cdot
 \widetilde{\alpha}(x_1) - (\mu)\cdot \widetilde{\alpha}(x_1) \\
 &= (\mu \phi_1)\cdot \tilde{R}(f_1,\dots,f_t) +
 \tilde{R}(\mu\phi_1,f_1,\dots,f_t) - (\mu)\cdot \tilde{R}(f_1,\dots,f_t)\\
 &= -(f_1)\cdots (f_t) + [-1/\mu]\cdot \tilde{R}(f_1,\dots,f_t) \\
 &\in \text{image}(\widetilde{\alpha}). 
}   
\end{small}
\end{proof}
\end{thm}

We know that $H^p(K,\Z(q))=0$ if $p>q$,
and we have the isomorphism (\ref{SVM}) with Milnor $K$-theory.
Therefore theorem (\ref{thm}) implies the following corollary.

\begin{corollary}\label{cor:ppTerme}
The morphism of $K_{*}^M(K)$ algebras 
\nneq
{
 K_*^M(K)\{\G_m(U)\}/J_U \xr{} \oplus_p \Hm{p}{U}{p} ; \quad (f) \mapsto [f],
} 
is an isomorphism, where the ideal $J_U$ is generated by elements of the form (\ref{Gl0}), (\ref{Gl1}) and (\ref{Gl2}).
\end{corollary}

\subsection{Relation with topological cohomology}

In the case $K=\C$ we have a realization functor $DM_{gm, \Q}^{eff,op}\xr{} D^{+}(\text{Vec}_{\Q})$ to the 
derived category of $\Q$ vector spaces (\cite{Huber}, 2.1.7),
which maps $M_{gm}(X)$ to the singular cochain  complex of $X^{an}$ and $\Q(1)$ to $H^1(\G^{an}_m,\Q)\cong \Q$. 
In particular we get a morphism of rings  
\eq{Real}
{
 \oplus_{p} H^p(U,\Q(p)) \xr{} \oplus_{p} H^p(U^{an},\Q).
}
We denote by $S(U)$ the ring $\oplus_{p} H^p(U,\Q(p))$. 
\begin{proposition}\label{prop:topcoh}
The homomorphism (\ref{Real}) induces an isomorphism of rings
\small
\nneq
{
 S(U)/\left(H^1(K,\Q(1)) \cdot S(U) \right) \xr{\cong}  \oplus_{p} H^p(U^{an},\Q).
} \normalsize
\begin{proof}
Using Proposition (\ref{prop:reinpp}) there is an isomorphism $ M_{gm}(U)\cong \oplus_{i\in I} \Z(n_i)[n_i]. $
This implies that $S(U)\cong \oplus_{i\in I} S(K)[-n_i] $ as $S(K)$ modules and $H^*(U^{an},\Q)\cong \oplus_{i\in I} \Q[-n_i]$, 
the map (\ref{Real}) beeing compatible with the decomposition by functoriality. Thus it is sufficient to prove 
$S(K)/\left(H^1(K,\Q(1)) \cdot S(K) \right)\cong \Q$, which follows from the comparison isomorphism (\ref{SVM}) since
$K^M_*(K)$ is generated by $K^M_1(K)$ as algebra.
\end{proof}
\end{proposition}


Using the realization functor of Ivorra \cite{Ivorra}: 
$$DM_{gm}^{eff}(K)^{op} \xr{} D^{+}(K,\Z_l)$$ 
to the triangulated category of $l$-adic sheaves defined by Ekedahl \cite{Ekedahl},  the statement of Proposition (\ref{prop:topcoh}) 
for $l$-adic cohomology can be proved  in the same way, with the assumption that $K$ is an algebraically closed field of characteristic 
$p\not= l$.

\subsection{Combinatorial description}
 When the ground field $K$ is the field of complex numbers $\C$, we have a combinatorial description of singular cohomology 
 at disposal (\cite{OrlikSolomon}, Theorem 5.2). 
 In this section we prove that this description holds for the ring  
 \nneq
 {
 A_0:=  \oplus_{p} H^p(U,\Z(p)) / \left(H^1(K,\Z(1)) \cdot \oplus_{p} H^p(U,\Z(p)) \right)       
 }
 for every perfect ground field $K$ as well. 

Let $Q$ be the cokernel of the inclusion $\G_m(K) \subset \G_m(U)$; by taking
divisors $\G_m(U)\xr{\text{div}} \oplus_{i=1}^r \Z \cdot Y_i$ in $\A^N$, we
get an isomorphism \\ $Q\xr{\cong} \oplus_{i=1}^r \Z \cdot Y_i$. We denote by 
$\Lambda_{\Z} Q$ the exterior algebra. Let $L$ be the ideal in $\Lambda_{\Z} Q$
generated by the elements 
\begin{small}
\nnal
{
 Y_{i_1}\wedge \dots \wedge Y_{i_s}; &\quad \text{if $Y_{i_1}\cap \dots \cap
                                                                 Y_{i_s}=\emptyset$} 
\\
 \sum_{k=1}^s (-1)^k Y_{i_1}\wedge \dots \wedge \widehat{Y_{i_k}} \wedge \dots \wedge Y_{i_s}; &\quad
 \text{if $Y_{i_1}\cap \dots \cap Y_{i_s}\not= \emptyset$} 
\\
 &\quad \text{and $\codim \left(Y_{i_1}\cap \dots \cap Y_{i_s}\right)<s$.} 
}
\end{small}
  
\begin{proposition}
Let $K$ be a perfect field. The map 
\nneq
{
\psi: A_0 \xr{} \Lambda_{\Z} Q /L ; \quad [f] \mapsto {\rm div}(f)
}
is well-defined and an isomorphism of rings.
\begin{proof}
By Corollary (\ref{cor:ppTerme}) we have
\nneq
{
 A_0 = K_*^M(K)\{\G_m(U)\}/\bigl(J_U+ K^{\times}\cdot K_*^M(K)\{\G_m(U)\}\bigr).
}
By definition, $\psi$ induces the projection $K_*^M(K)\xr{} K^M_0$ 
on $K_*^M(K)$. We have to prove that elements of the form 
(\ref{Gl1}) and (\ref{Gl2}) map to zero. For elements of the
form (\ref{Gl2}) the assertion is trivial and we may replace (\ref{Gl1})
with (\ref{Gl2R}) and we may assume $f_j=\lambda_j\phi_{i_j}$ or $f_j=\lambda_j$  for
every $j$ and some $\lambda_j\in K^{\times}$, as in the proof of the theorem (\ref{thm}).
We need to prove that 
\nneq
{
\alpha:= \sum_{k=1}^s (-1)^k {\rm div}(f_1)\wedge \dots \wedge \widehat{{\rm div}(f_k)} \wedge 
                     \dots \wedge {\rm div}(f_s) 
}
is an element of $L$. Obviously, $\alpha$ is trivial if there are more than two constant functions among the
$f_1,\dots,f_s$; if $f_1$ is the only constant function then $Y_{i_2}\cap
\dots \cap Y_{i_s}=\emptyset$ and $\alpha\in L$. In the case of no non-constant
function, then either $Y_{i_1}\cap \dots \cap Y_{i_s}=\emptyset$ or 
$\codim \left(Y_{i_1}\cap \dots \cap Y_{i_s}\right)<s$. In the first case we
have $Y_{i_1}\cap \dots \cap \widehat{Y_{i_j}}\cap \dots \cap
Y_{i_s}=\emptyset$ for every $j$ and in the second case $\alpha\in L$ by
definition of $L$. This proves that $\psi$ is well-defined.

For the inverse map we have to prove that
\nneq
{ 
 \Lambda_{\Z} Q /L \xr{} A_0; \quad Y_i\mapsto [\phi_i],
}   
is well-defined. 
Because $A_0$ is graded-commutative, the map $\Lambda_{\Z} Q\xr{} A_0, Y_i\mapsto [\phi_i],$ exists. 
If $Y_{i_1}\cap \dots \cap Y_{i_s}=\emptyset$ then
$\sum_{k}\lambda_k \phi_{i_k}=1$ for some $\lambda_k\in K$ and we 
have $[\phi_{i_1}]\cdots [\phi_{i_s}]=0$ in $A_0$. If 
$Y_{i_1}\cap \dots \cap Y_{i_s}\not= \emptyset$ and $\codim \left(Y_{i_1}\cap \dots \cap Y_{i_s}\right)<s$
then $\sum_{k}\lambda_{k} \phi_{i_k}=0$ for some $\lambda_k\in K$ and 
$\sum_{k} (-1)^k [\phi_{i_1}]\cdots \widehat{[\phi_{i_k}]} \cdots
[\phi_{i_s}]=0$ in $A_0$. This proves that the inverse map is well-defined. 
\end{proof}
\end{proposition}

\begin{corollary}\label{Rang}
The rank of the free $H(K)$ module $H(U)$ is equal to the rank of $\Lambda_{\Z} Q/L$.
\begin{proof}
The rank of $H(U)$ is equal to the rank of $\oplus_{p} H^p(U,\Z(p))$ as 
$\oplus_{p} H^p(K,\Z(p))\cong K^M_{*}(K)$ module, by Proposition (\ref{prop:reinpp}). 
Thus, we get
\nneq
{
 \rk_{K^M_{*}(K)}(\oplus_{p} H^p(U,\Z(p))) = \rk_{\Z} A_0 = \rk_{\Z} \left(\Lambda_{\Z} Q/L\right).
}
\end{proof}
\end{corollary}

The algebra $\Lambda_{\Z} Q/ L$ depends only on the combinatorics of the hyperplanes
in the complement. For example, assume that
\eq{BedNKD}
{
 \overline{Y_1}, \dots ,\overline{Y_r}, H \quad \text{is a normal crossing divisor,}
}
where $H$ is the hyperplane at infinity for the compactification $\A^N\subset
\P^N$ and $\overline{Y_i}$ is the closure of $Y_i$ in $\P^N$. Then we have
\nneq
{
 \Lambda_{\Z} Q/L \cong \Lambda_{\Z} Q / \langle Y_{i_1}\wedge \dots \wedge
 Y_{i_{N+1}} \mid i_1,\dots,i_{N+1}\rangle,
}
here $\Lambda_{\Z} Q/L$ depends only on $(r,N)$.

The motivic cohomology ring is a finer invariant. It is easy to see (with
assumption (\ref{BedNKD})) that 
\eq{Basis}
{
\{[\phi_{i_1}]\cdots [\phi_{i_t}] \mid i_1<\dots<i_t, t\leq N \} \cup \{1\}
} 
is a base of $ H(U)$ as module over the motivic cohomology of the ground field. We can define an homomorphism
of groups 
\eq{tames}
{
 \beta: \G_m(U)^{\otimes N+1} \xr{} K^{\times} \otimes_{\Z} \Lambda^N Q
} 
by setting 
\small 
\nneq
{
\beta(f_1\otimes \cdots \otimes f_{N+1})=\sum_{i_1<\dots<i_{N}}
\alpha_{i_1<\dots <i_{N}}\otimes \left(Y_{i_1}\wedge \dots \wedge Y_{i_N}\right) 
}
\normalsize
with 
\nneq
{
 [f_1]\cdots [f_{N+1}]= \sum_{i_1<\dots<i_{N}}
    [\alpha_{i_1<\dots <i_{N}}] \cdot [\phi_{i_1}]\cdots [\phi_{i_N}] + \dots 
}
 where this expression is the presentation of $[f_1]\cdots [f_{N+1}]$ in the base (\ref{Basis}). The
morphism $\beta$ is independent of the choices for $\phi_1,\dots,\phi_{r}$. 

For $U=\A^1-\{p_1,\dots,p_r\}$ it is easy to see that $\beta$ is the ``tame
symbol'':
\nneq
{
 \beta(f,g)= \sum_{i=1}^r (-1)^{\nu_{p_i}(f)\nu_{p_i}(g)} \left(
 \frac{f^{\nu_{p_i}(g)}}{g^{\nu_{p_i}(f)}} \right)(p_i) \otimes p_i.
}

\end{document}